\documentclass[a4paper,12pt]{article}

\usepackage{amsmath}
\usepackage{amsthm}
\usepackage{amsfonts}
\usepackage{amssymb}

\numberwithin{equation}{section}

\newtheorem{theorem}{Theorem}[section]
\newtheorem{cor}[theorem]{Corollary}
\newtheorem{lemma}[theorem]{Lemma}
\newtheorem{prop}[theorem]{Proposition}

\theoremstyle{definition}

\newtheorem{defi}{Definition}[section]
\newtheorem{rem}[defi]{Remark}

\newtheorem{vuoto}[defi]{Fact}

\newcommand{\ov}{\overline}

\newcommand\lan{\langle}
\newcommand\ran{\rangle}
\newcommand\st{{\lan\cdot,\cdot\ran}}

\newcommand{\half}{{\frac{1}{2}}}
\newcommand\psuc{{\frac{p}{c}}}

\newcommand{\Q}{Q^\vee}
\newcommand\Qst{\Q_\st}
\newcommand{\Dp}{\Delta^+}
\newcommand{\D}{\Delta} 
\newcommand{\Wa}{\widehat{W}}
\newcommand{\wa}{\widehat w}
\newcommand{\ua}{\widehat u}

\newcommand{\ww}{{w_*}}
\newcommand{\iw}{i^{w_*}}

\renewcommand{\th}{\theta}

\renewcommand{\l}{\lambda} 
\renewcommand{\b}{\beta}

\renewcommand{\r}{\rho} 
\renewcommand{\a}{\alpha} 
 
\renewcommand{\o}{\omega}

\newcommand{\s}{\sigma}

\newcommand{\h}{\mathfrak{h}}
\newcommand{\g}{\mathfrak{g}}

\newcommand{\R}{\mathbb{R}}
\newcommand{\real}{{\mathbb{R}}}
\newcommand{\nat}{{\mathbb{N}}}
\newcommand{\ganz}{{\mathbb{Z}}}

\date{}

\begin{document}

\title
{The $\Wa$-orbit of $\r$, Kostant's formula for powers of
the Euler product and affine Weyl groups as permutations of $\ganz$}

\author{Paola Cellini\\ Pierluigi M\"oseneder Frajria\\Paolo Papi}

\maketitle

\begin{abstract} 
Let an affine Weyl group $\Wa$ act as a group of affine transformations on a 
real vector space $V$.
We analyze the $\Wa$-orbit of a regular element in $V$ and deduce
applications to Kostant's formula for powers of
the Euler product and to the  representations of $\Wa$ as permutations of the 
integers.
\end{abstract}

\section{Introduction}

This paper stems from the attempt of deepening  two seemingly unrelated topics:
on  one hand the combinatorial interpretation of  Kostant's recent results on 
the powers of the
Euler product suggested in type $A$  by Tate and Zelditch, on the other hand the 
problem of giving
a uniform and conceptual description of  certain affine Weyl groups  as
permutations of the  integers. The common denominator of these two  subjects 
lies in their
intimate connection with the orbit  of a
distinguished vector under the action of an  affine Weyl group.
The results of the paper should be regarded for the first topic as a generalization
 of Tate and Zelditch's approach, for the other
as a systematic treatment of well-established results on  affine Weyl groups of classical type. 
To be more precise let us fix
notation. Let
$(V,(\cdot, \cdot))$ be an Euclidean space,
$\D$ a finite crystallographic  irreducible root system in $V$, $\Dp$ a fixed 
positive system for
$\D$.\par Set $\rho=\half\sum\limits_{\a\in\Dp}\a$ and  let $\th$ be the highest 
root of $\D$.
We define the dual Coxeter 
number $h^\vee$ of $\D$ as 
$
h^\vee=\frac{2(\rho,\theta)}{(\theta,\theta)}+1.
$
The affine Weyl group $\Wa$ of $\D$ is the group 
generated by reflections on $V$ with respect to the set of affine hyperplanes 
$H_{\a,k}=\{x\in V\mid (x,\a)=k\}$, $\a\in \Dp,\,k\in\ganz$. For each $q\in 
\real^+$, we denote by $\Wa_q$ the group generated by reflections in $V$ with 
respect to the set of hyperplanes $H_{\a,qk}$, $\a\in \Dp,\,k\in\ganz$; 
thus $\Wa_q$ is naturally isomorphic to $\Wa$. 
We notice that scaling the inner product by $\frac{1}{q}$ changes $\Wa$ into 
$\Wa_q$ (and does not change $h^\vee)$. 
We assume along the paper that 
\begin{equation}\label{norm}
(\th, \th)=\frac{1}{h^\vee}.
\end{equation}
For each $\l\in V$, we denote by $\Wa_q\cdot \l$ the orbit of $\l$ under 
$\Wa_q$.
\par
A basic step of our work is the analysis of $\Wa_\half\cdot \r$.
A motivation for this study occurs in the framework of Kostant's work 
on Dedekind's $\eta$ function, which we now recall. Let $\g$ be a complex finite 
dimensional 
semisimple Lie algebra, $\h$ a Cartan subalgebra of $\g$ and $\D$ the 
corresponding root system. Let $V=\h_\R^*$, the real span of a fixed set of 
simple roots, endowed with the invariant form induced by the Killing form of 
$\g$. (It is well-known that then \eqref{norm} holds).
\par 
If $\l$ is a dominant weight let $\chi_\l$ denote the character of the 
irreducible $\g$-module $V_\l$ with highest weight $\l$. Set also 
$a=\exp(2\pi i\cdot 2\rho)$. Working on previous results of Macdonald, 
Kostant found the following remarkable expansion for (certain) powers of the 
Euler product $\prod_{m=1}^\infty(1-x^m)$. 
\begin{theorem}
\cite[Thm 3.1]{kostantadv}
\begin{equation}
\left(\prod_{m=1}^\infty(1-x^m)\right)^{\dim(\g)}=\sum_{\l\text{
dominant}}\chi_\l(a)\,\dim(V_\l)\,x^{(\l+2\rho,\l)}.
\label{formula}
\end{equation} 
Moreover, $\chi_\l(a)\in\{-1,0,1\}$.
\end{theorem}

\vspace{5pt}

In \cite{kostant} Kostant has improved the previous formula determining the 
set $P_{alc}$ of weights which give non zero contribution in the sum (see 
Theorem \ref{kost} below). The main outcome is that
$$
P_{alc}=\{\l\text{ dominant weight}\mid \l+\r\in\Wa_\half\cdot\rho\}.
$$
Moreover, he proves that the contribution of each $\l\in P_{alc}$ is 
determined by the parity of $\ell_\half(w)$, where $w\in \Wa_\half$ is the 
element 
such that $\l+\r=w(\r)$, and $\ell_\half$ is the length function on $\Wa_\half$.    

\vspace{5pt}
On the other hand, in \cite{adin}, Adin and Frumkin  made 
explicit, by using the well-known connection between dominant weights and 
partitions, the combinatorial content of Kostant's result in type $A$. Their 
result makes also easy to determine the sign of $\chi_\l(a)$. After the 
appearance of Kostant's paper, a simple approach to the combinatorial 
interpretation of Kostant's result in type $A$ using the affine Weyl group was 
explained by Tate and Zelditch in \cite{zelditch}.
We shall obtain results analogous to those of \cite{zelditch} for 
all classical types and for $G_2$. The exposition of these results is 
the content of Section \ref{esempi}. 

\par
The crucial observation is that $\r$ is the {\it unique} weight of 
$\g$ lying in the fundamental alcove of $\Wa_\half$. By the basic properties 
of the action of the affine group on $V$, this implies that 
$\Wa_\half\cdot\r$
is the set of weights which lie in some alcove of $\Wa_\half$, or, equivalently, 
which do not belong to any of the reflecting hyperplanes. Once the root systems
are explicitly described in coordinates, this allows us to easily describe 
$P_{alc}$ by purely arithmetical conditions, for all types. 
\par
We shall write down this description only for the classical types and for 
$G_2$. For each of these cases, we shall also give a simple rule for 
recovering the parity of $\ell(w)$ from $w(\r)$. For type $A$, we re-obtain 
the rule of \cite{adin}. The affine Weyl group is the semidirect product 
of the finite Weyl group $W$ of $\g$ and the group $\Q$ acting on $V$ by 
translations, hence $\Wa_\half\cong \half\Q \rtimes W$. Moreover, if 
$w=t_\tau v$,where $t_\tau$ is the translation by $\tau\in \half\Q$, and 
$v\in W$, then $\ell(w)\equiv \ell(v)\mod 2$. Our rule is in fact a sort of 
Euclidean algorithm which produces $v$ and $\tau$ from $w(\r)$. 

\vspace{5pt}

The last section of the paper deals with affine Weyl groups regarded as 
permutation groups of the set of integers. This point of view has 
been introduced by Lusztig \cite{L} for type $\tilde A$, and generalized to 
the other classical cases by his students (and other people). A thorough 
and systematic account of the combinatorial aspects of the theory can now be 
found in Chapter 8 of  \cite{BBlibro}. 
\par

From the explicit description of $\Wa_\half \cdot \rho$, we see that in
cases 
$\tilde A$ and $\tilde C$ we can quite naturally associate to each $w\in 
\Wa_\half$ a permutation of $\ganz$, uniquely determined by $w(\r)$. In 
this way, we obtain an injective homomorphism of $\Wa_\half$ into $S(\ganz)$, 
the group of permutations of $\ganz$, which agrees with the usual permutation 
representation. This suggests that the known permutation representations of 
all classical affine Weyl groups can be obtained from the explicit 
description the orbit $\Wa_q\cdot \l$, for an appropriate choice of $q$ and 
$\l$.
In fact, the final outcome of our study is a uniform and concise treatment of 
the known permutation representations of classical Weyl groups. Our point of 
view is also successful for type $\tilde G_2$. To our knowledge, a similar 
unified approach does not appear in literature, even if the existence of a 
connection between the orbit of a regular vector and the permutation 
representation of $\Wa$ is noted in H.~Eriksson's unpublished Ph.D. Thesis 
\cite{Eth}.\par
We have already explained the content of Sections 3, 4. The results of Section 2 
are 
a kind of  ``context free" preparation to the next Sections, and rely on the 
standard theory 
of the geometric action of  affine Weyl groups. The main contribution
is Proposition \ref{fund}.

\section{Preliminary results.}
\label{orbita}

We retain the notation set at  the beginning of the Introduction:  $V$ is an 
$n$-dimensional Euclidean space with inner product $(\cdot,\cdot)$,  $\D$ is a 
finite crystallographic irreducible root system of rank $n$ in $V$. Denote by 
$W$  the corresponding finite reflection group. Let $\Pi=\{\a_1,\ldots,\a_n\}$ 
be a set of simple roots for $\D$ (with  positive system $\Dp$).
 Denote by
$Q$ the root lattice. For $\b\in Q$ set $\b^\vee=\frac{2\b}{(\b,\b)}$,  and let
$$\begin{array}{lcl}
Q^\vee & = &\sum\limits_{i=1}^n\ganz \a_i^\vee,
\\ \\
P & = &\left\{\l\in\h\mid
(\l,\a^\vee)\in\ganz\,\forall\,\a\in\D\right\},
\end{array}$$
be the coroot and weight lattices. Denote by  $P^+$  the 
set of dominant weights:
$$P^+=\{\l\in P\mid (\l,\a^\vee)\geq 0\ \forall\,\a\in\Pi\},
$$
Let  $\omega_1,\ldots,\omega_n$ be the fundamental weights, so that 
$P=\sum\limits_{i=1}^n\ganz
\omega_i$ and $\rho=\sum\limits_{i=1}^n
\omega_i$. Remark that if  
$\th^\vee=\sum\limits_{i=1}^nm_i\a_i^\vee$ then $h^\vee= 
1+\sum\limits_{i=1}^nm_i$.

Fix $q\in\R^+$. Recall the group $\Wa_q$ defined in the Introduction. Then $\Wa_q=T_q\rtimes W$ 
where
$T_q$ is the group of translations of $V$ by elements in $qQ^\vee$. It is clear that $\Wa_1$ is the usual affine Weyl group.
Ours is a slight extension of the usual definition 
which turns out to be very useful for our goals.
\par 
For $\a\in V\setminus\{0\},\,\b\in V$ denote by $s_\a,\,t_\b$ the reflection 
in $\a$ and the translation by $\b$, respectively.
\par
Recall that $\Wa_q$ is a Coxeter group with generators
$s_i=s_{\a_i}$ for $i=1,\dots,n$ and $s_0=t_{q\th^\vee}s_\th$. We
denote by $\ell_q$ the length function with respect to this choice
of generators. Set $H_{rq,\a}=\{x\in V\mid (x,\a)=rq\}$ for
$r\in \ganz$ and $\a\in\Dp$. The {\it alcoves} of $\Wa_q$ are 
the connected components of  $V\setminus
\bigcup_{\substack{\a\in\Dp \\ r
\in\ganz}}H_{rq,\a}$. 
The {\it fundamental alcove} is the alcove  
$$
C_q=\left\{x\in V\mid (x,\a)> 0\,\forall\,\a\in\Dp,\,\,(x,\th)<q\right\}.
$$
It is well-known that $\Wa_q$ acts on the set of alcoves and this action is 
simply transitive. 
This means that $wC_q$ is an alcove and, for each alcove $C_q'$ there exists a
unique
$w\in
\Wa_q$  such that $C'_q=w(C_q)$. Moreover, 
$\overline{C_q}$ is a fundamental domain for the action of
$\Wa_q$ on $V$. In particular if $y$ belongs to some alcove, then there exist unique
$w\in\Wa_q$ and $x\in C_q$ such that $w(x)=y$. We shall  tacitly use these standard properties in the
following.

\begin{defi} 
We say that $v\in V$ is $q$-regular if it belongs to some alcove, 
or, equivalently,
$$
v\in V\setminus\bigcup_{\substack{\a\in\Dp \\ r\in\ganz}}H_{rq,\a}.
$$
\end{defi} 

Any alcove can be expressed as an intersection (ranging over $\Dp$)
of strips
$H_\a^{rq}=\{x\in V\mid rq<(x,\a)<(r+1)q\},$
$(r\in\ganz)$. Denote by   $k(w,\a)$ the integers such that
$$
wC_q=\bigcap_{\a\in\Dp} H_\a^{k(w,\a)q}.
$$
The collection $\{k(w,\a)\}_{\a\in\Dp}$ has been introduced by Shi and
called the {\it alcove form} of $w$.

\begin{rem}
\label{regolare} 
Suppose that $\mu$ is $q$-regular. If $\mu\in wC_q$, then

\begin{equation}
\label{alcove}
k(w,\a)=\left\lfloor \frac{(\mu,\a)}{q}\right\rfloor
\end{equation}
and

\begin{equation}\label{lunghezza}
\ell_q(w)=\sum_{\a\in\Dp}|\left\lfloor 
\frac{(\mu,\a)}{q}\right\rfloor|.
\end{equation}

To obtain \eqref{alcove}, it suffices to remark that the r.h.s. counts the
number of hyperplanes   $H_{rq,\a}$ separating $C_q$ from
$wC_q$. Since the total number of separating hyperplanes  $H_{rq,\a}$ when $\a$ 
ranges over $\Dp$, gives $\ell_q(w)$ (see \cite[4.5]{H}), \eqref{lunghezza} 
follows.
\end{rem}
\vskip5pt

We state as a proposition the following elementary observation, which will 
play a prominent role in the sequel.

\begin{prop}
\label{fund} 
Fix $\l\in V$. Let $L$ be a lattice in $V$ such that $\l+L$ is 
$\Wa_q$-stable and $(\l+L)\cap C_q=\{\l\}$. Then 
$$\Wa_q\cdot\l= \{\mu\in\l+L\mid \text {for all } \a\in\D, 
\ (\mu, \a)\notin q\ganz\}.
$$
\end{prop}

\begin{proof} 
Assume $\mu\in \Wa_q\cdot\l$. Since $\l+L$ is $\Wa_q$-stable, $\mu\in\l+L$. 
Moreover, since $\Wa_q$ acts on the set of alcoves,    
$\mu$ belongs to some alcove, which means that for all $\a\in\D$, 
we have $(\mu, \a)\notin q\ganz$.
\par 
Conversely, assume that $\mu\in\l+L$ and, for all 
$\a\in\D, \ (\mu, \a)\notin q\ganz$. Then $\mu$ belongs to some alcove. 
Since $\Wa_q$ acts transitively on the set of alcoves, and preserves
$\l+L$, there exists $w\in \Wa_q$ such that $w(\mu)\in C_q\cap (\l+L)=\{\l\}$.
\end{proof}

\begin{rem} 
If $\l+L$ is $W$-stable and $qQ^\vee\subset L$ then $\l+L$ is $\Wa_q$-stable.
\end{rem}

\begin{lemma}
\label{taut} 
We have $C_q\cap P=\{\r\}$ if and only if 
\begin{equation}
\label{eqtaut}
\frac{(\th,\th)}{2}(h^\vee-1)<q\leq \frac{(\th,\th)}{2}(h^\vee+m-1)
\end{equation} 
where $m=\min_{1\leq i\leq n} m_i$. In particular, 
$$
P\cap C_\half=\{\rho\}.
$$
\end{lemma}

\begin{proof} 
Note that $(\rho,\th)=\frac{(\th,\th)}{2}(h^\vee-1)$, hence 
$\rho\in C_q\cap P$ if and only if $\frac{(\th,\th)}{2}(h^\vee-1)<q$.
Obviously $C_q\cap P=\{\rho\}$ if and only if $\rho+\omega_i\notin C_q$ for all 
$i=1,\ldots,n$.
This implies
$$q\leq (\rho+\omega_i,\theta)=\frac{(\th,\th)}{2}(h^\vee-1)
+\frac{(\th,\th)}{2}m_i=\frac{(\th,\th)}{2}(h^\vee+m_i-1)$$
as desired.
\end{proof}
 Note that $m=1$ if $\D$ is not of type $E_8$; in this latter case 
$m=2$.

\section
{Application to Euler products.}
\label{esempi}

The first application of the above results is connected with the work of Kostant 
on the powers of
the Euler product
$\prod_{m=1}^\infty(1-x^m)$.
\par
Let $\g$ be a complex finite dimensional semisimple Lie algebra, $\h$  a Cartan
subalgebra of $\g$ and $\D$ the corresponding root system. In the notation of 
the previous
section, we choose
$V$ to be the real span $\h_\R^*$ of a fixed set of simple roots endowed with 
the invariant 
form induced
by the Killing form of $\g$. With this choice we have indeed that
$(\th,\th)=\frac{1}{h^\vee}$ (see e.g. \cite[\S~2]{kostantadv}).
\par If
$\l\in P^+$, let $\chi_\l$ denote the  character of the
irreducible $\g$-module $V_\l$ with highest weight $\l$. 
Recall relation \eqref{formula}. In \cite[Theorem 2.4]{kostant} 
a general criterion for determining the set
$$P_{alc}=\{\l\in P^+\mid \chi_\l(a)\ne0\}$$
is provided
 (see also \cite[Exercise 10.19]{Kac}).
Kostant's theorem  can be rephrased as follows:

\begin{theorem}\label{kost} We have
$$\l\in P_{alc}\iff  \l+\r\in
\Wa_{\half}\cdot\rho.$$ Moreover, if
$\l+\r=w(\rho),\,w\in\Wa_{\half},$ then
$\chi_\l(a)=(-1)^{\ell_{\frac{1}{2}}(w)}.$
\end{theorem}

\begin{cor}\label{perkostant} A weight  $\l$ belongs to $P_{alc}$ if
and only if it is dominant and
\begin{equation}\label{cond}
(\l+\r,\a)\notin
\frac{1}{2}\ganz\quad\text{for any $\a\in\D$.}
\end{equation}
In such a case, $\l$ belongs to the root lattice $Q$ and
\begin{equation}\label{segno}
\chi_\l(a)=(-1)^{\sum_{\a\in\Dp}\lfloor
2(\l+\r,\a)\rfloor}.
\end{equation}
\end{cor}
\begin{proof}
By Lemma \ref{taut} we have that  $C_{\half}\cap P=\{\rho\}$.
Recall that $(\th,\th)=\frac{1}{h^\vee}$. 
Then
$\frac{1}{2}Q^\vee\subset Q$, hence we can apply
Proposition \ref{fund}.
Moreover, if $\l+\r\in \Wa_{\tfrac{1}{2}}\cdot\rho$, then
$\l+\r\in \r+Q+\frac{1}{2}Q^\vee\subset \r+Q$, hence $\l\in
Q$.
Finally  \eqref{segno} follows readily from Theorem \ref{kost} and
\eqref{lunghezza}.
\end{proof}

In the rest of this section we provide an explicit rendering of
Corollary \ref{perkostant}  for the classical root systems.  We find
combinatorial conditions  that guarantee that $\l\in P_{alc}$ and
determine $\chi_\l(a)$. For this last purpose is convenient to use the following 
general fact
rather than formula \eqref{segno}. Denote by $\ell$ the length function in $W$.
\begin{lemma}\label{lunghezzatau}  If $t_\tau w\in\Wa_q,\,\tau\in q\Q,w\in W$, 
then 
$\ell_q(t_\tau
w)\equiv\ell(w)\,\mod\,2$.
\end{lemma}
\begin{proof} We shall use several times the following well-know fact from the 
theory of Coxeter
groups (see e.g. \cite[5.8]{H}): cancellations occur in pairs, so that if an 
element has an expression 
in terms of the generators of a certain parity, its length has the same parity.
Since $t_\tau w$ has certainly an expression involving $\ell_q(t_\tau)+ \ell(w)$ 
generators, it suffices 
to show that $\ell_q(t_\tau)$ is even. Since $q\Q$ is the $\ganz$-span of 
$qW\cdot\theta^\vee$
it suffices to prove that if $u\in W$, then $\ell_q(t_{qu(\theta^\vee)})$ is 
even.  This   follows from the relation
$t_{qu(\theta^\vee)}=us_0s_\theta u^{-1}$.
\end{proof}

In the classical cases we shall explicitly determine for each $\l\in P_{alc}$ 
the 
unique element $w\in
\Wa_{\frac{1}{2}}$ such that $\l+\rho=w\rho$ and compute $\tau\in 
\frac{1}{2}\Q$, $u\in W$ such that $w=t_\tau u$.
Applying Lemma
\ref{lunghezzatau} we obtain that $\chi_\lambda(a)=(-1)^{\ell(u)}$. 
 In
\cite{zelditch} essentially
 the same analysis was applied only to type $A_{n}$ obtaining Theorem 1.2 of 
\cite{adin}. 
In the following we adopt the realization of the irreducible root
systems as subsets of $\R^N$ given in \cite{bourbaki}. We denote
by $\langle\cdot,\cdot\rangle$ the standard inner product of $\R^N$ and by
$\{e_i\}$ the canonical basis.

   \subsection{Type $A_{n}$.}\label{tipoA}
  Recall
  that in  \cite{bourbaki}
  $\h^*_\R$ is identified with the subspace of $\R^{n+1}$
  orthogonal to $\l_0=\sum\limits_{i=1}^{n+1}e_i$. In this setting
  $$
\Dp=\{e_i-e_j\mid i<j\}
  $$
  and
  $$
  Q=(\sum\limits_{i=1}^{n+1}\ganz e_i)\cap \h_\R^*.
  $$
  The map  $\l\mapsto\bar{\l}= \l-\langle\l,e_{n+1}\rangle\l_0$  maps
  $P$ bijectively onto
$\sum\limits_{i=1}^n\ganz e_i$, $P^+$ onto
$$
P_n=\left\{\sum_{i=1}^n \l_i e_i \mid \l_i\in\ganz,\
\l_1\ge\l_2\ge\dots\ge\l_n\ge0\right\}.
$$
We finally recall that  $\r=\sum\limits_{i=1}^{n+1}
\frac{n-2i+2}{2}e_i,\,\th=e_1-e_{n+1}$, hence $h^\vee=n+1$ .
Since $\langle\th,\th\rangle=2$ and	$(\th,\th)=\frac{1}{h^\vee}$, we have  
\begin{equation}\label{killinga}(\cdot,\cdot)=\frac{1}{2h^\vee}
\langle\cdot,\cdot\rangle.\end{equation}
This implies in particular that $\frac{1}{2}\Q=(n+1)Q$.\vskip5pt
If $\l\in\h^*_\R$ set $\l_i=\langle\bar{\l},e_i\rangle$.
  Since $\langle\l_0,\a\rangle=0$ for all $\a\in \h^*_\R$
we see that $\langle\bar{\l},\a\rangle=\langle\l,\a\rangle$ for all $\a\in \D$.
Also recall that 
$\bar{\r}=\sum\limits_{i=1}^n (n-i+1)e_i$. Applying Corollary \ref{perkostant} 
we deduce the following result,
which is the first statement of Theorem 1.2 from
\cite{adin}.
\begin{prop}\label{tipoa} For $\ov\l=\sum\limits_{i=1}^n\l_i e_i\in P_n$ we have
\begin{equation*} 
\l\in P_{alc}\iff
 \l_i+n-i+1\ne \l_j+n-j+1\
mod\,(n+1).
\end{equation*}
$(1\leq i\ne j\leq n+1)$.\end{prop} 
\vskip10pt
Note that, since $\l\in Q$, we have  
$$\sum_{i=1}^{n+1}\l_i=\sum_{i=1}^{n+1}\langle \ov\l,e_i\rangle=
(\sum_{i=1}^{n+1}\langle \l,e_i\rangle)-(n+1)\langle \l,e_{n+1}\rangle=-
(n+1)\langle \l,e_{n+1}\rangle.$$ 
 Hence $n+1$
divides
$\sum\limits_{i=1}^{n+1}\l_{i},$ so we can
write
\begin{equation}\label{divisioneA}
\l_{i}+(n-i+1)-\frac{1}{n+1}\sum_{j=1}^{n+1}\l_{j}=(n-r_{i}+1)+(n+1)q_{i}
\end{equation} with
$r_{i}\in\{1,2,\dots,n+1\}$.

Set $\tau=(n+1)\sum\limits_{i=1}^{n+1}q_{i}e_{i}.$
By Proposition \ref{tipoa} the $r_{i}$ are pairwise distinct, so,
by \eqref{divisioneA}
$$
(n+1)\sum_{i=1}^{n+1}q_{i}=\sum_{i=1}^{n+1}(n-i+1)-\sum_{i=1}^{n+1}(n-
r_{i}+1)=0,
$$ hence $\tau\in\frac{1}{2}\Q$.
We can write
\begin{align*}
\l+\r&=\sum\limits_{i=1}^{n+1}\left(\l_{i}+(n-i+1)-\frac{1}{n+1}\sum\limits_{j=1}^{n+1}\l_{j}-
\frac{n}{2}\right)e_{i}\\
&=\sum\limits_{i=1}^{n+1}(\frac{n-2r_{i}+2}{2})e_{i}+(n+1)
\sum\limits_{i=1}^{n+1}q_{i}e_{i}.
\end{align*}
The action of $W$ on $V$ is described explicitly in \cite{bourbaki}. In particular it is known that, if $v\in W$, then there is an element $\s_v$ of $S_n$ such that $v(e_i)=e_{\s_v(i)}$. This fact establishes the well known isomorphism between $W$ and $S_n$. Thus if we set $\s$ to be the element of $S_n$ such that $\s(i)=r_i$, and let $v$ be the 
element of $W$ such that $\s_v=\s^{-1}$,  then
$v(\r)=\sum\limits_{i=1}^{n+1}(\frac{n-2i+2}{2})v(e_{i})=
\sum\limits_{i=1}^{n+1}(\frac{n-2\s(i)+2}{2})e_{i}=
\sum\limits_{i=1}^{n+1}(\frac{n-2r_{i}+2}{2})e_{i}$
hence $\l+\r=t_\tau v(\r)$ and $\chi_\l(a)=(-1)^{\ell(v)}$.

\begin{rem}
It is well known (and easy to prove) that $(-1)^{\ell(v)}=sign(\s_v)$ thus $\chi_\l(a)$ is the sign of the permutation $i\mapsto r_i$.
\end{rem}

\subsection{Type $C_{n}$.}\label{tipoC}
We have $\Dp=\{e_i\pm e_j\mid 1\leq i<j\leq n\}\cup \{2e_i\mid 1\leq
i\leq n\},\,\r=\sum\limits_{i=1}^n (n-i+1)e_i,\,\th=2e_1$ so that  $h^\vee=n+1$. 
Moreover
\begin{align*}
&P=\sum_{i=1}^n\ganz e_i,\qquad Q=\{\sum_{i=1}^{n}\l_ie_i\mid
\sum_{i=1}^{n}\l_i\in 2\ganz\},
\\
&P^+=\left\{\sum_{i=1}^n\l_i e_i\in P\mid \l_1\geq\l_2\geq\dots\l_n\geq 
0\right\}.\end{align*}
This time $\langle\th,\th\rangle=4$, so that $(\cdot,\cdot)=\frac{1}{4h^\vee}
\langle\cdot,\cdot\rangle$ and
$\frac{1}{2}\Q=2h^\vee \ganz^n$. By Corollary \ref{perkostant} we have
  
\begin{prop}\label{tipoc} For $\l=\sum\limits_{i=1}^n\l_i e_i\in P^+$ we have 
$$\l\in P_{alc}\iff\begin{array}{l} \l_i+n-i+1\not\equiv \pm(\l_j+n-j+1)
\mod\,2(n+1)\quad\!\!(i\ne j)\\ \l_i+n-i+1\notin
(n+1)\ganz.
\end{array}$$
\end{prop}
\vskip5pt
It is well-known that the finite Weyl group $W$ acts faithfully on 
$\{\pm e_1,\ldots,$ $\pm e_n\}$ by signed permutations. It follows that $W\cdot 
\rho$ is the set
of elements of type $\sum\limits_{i=1}^n a_ie_i$ with $\{\pm a_1,\ldots,\pm 
a_n\}=
\{\pm 1,\ldots,\pm n\}$ . Now assume that $\l\in P_{alc}$ and $\mu=\l+\r$, 
$\mu=\sum\limits_{i=1}^n\mu_i e_i$. 
Denote by $\ov\mu_i$ the unique element in $\{\pm 1,\ldots,\pm n\}$
such that $\mu_i\equiv \ov\mu_i\,\mod 2(n+1)$ and set 
$\ov\mu=\sum\limits_{i=1}^n 
\ov\mu_ie_i$.
Notice that by Proposition
\ref{tipoc} the $\ov\mu_i$ are distinct and different from $0,n+1$.
Then there exists $v\in W$ such that $\ov\mu=v(\r)$. Moreover  from 
the description of $\frac{1}{2}\Q$ it follows that $\mu-\ov\mu\in\frac{1}{2}\Q$. 
Set $\tau=\mu-\ov\mu$.
It follows that $\l+\r=t_\tau v(\rho)$ and hence, by Lemma \ref{lunghezzatau}, 
we 
have
$\chi_\l(a)=(-1)^{\ell(v)}$. 

 \begin{rem}\label{segnoc} If $v\in W$ define   $(\pm i)^{\s_v}=\pm\left\langle v(\r),e_{n-i+1}\right\rangle$ for $i=1,\dots,n$. Since, as observed above,   $W$ acts as signed permutations on 
$\{\pm e_1,\ldots,\pm e_n\}$ we have that the map $v\mapsto \s_v$ defines an homomorphism from $W$ to the set of signed permutations on $\{\pm 1,\ldots, \pm n\}$. If $\s$ is such a signed permutation then set $|\s|$ to be the element of $S_n$ defined by $i^{|\s|}=|i^\s|$ and set $n_{\s}=\sharp\left\{i\mid i^{\s}<0,\, i=1,\dots,n\right\}$. It is well known that $\chi(\s)= sign(|\s|)(-1)^{n_\s}$ is a character of the group of signed permutations. Since $\chi(\s_{s_i})=-1$ it follows at once that $(-1)^{\ell(v)}=\chi(\s_v)$.
This shows that $\chi_\l(a)=sign(\left|\s_v\right|)(-1)^{n_{\s_v}}$.  Observe that $|\s_v|$ is the permutation of $\left\{1,2,\dots,n\right\}$ defined by setting $i^{\left|\s_v\right|}=|\ov\mu_{n-i+1}|$ and $n_{\s_v}=\sharp\left\{i\mid \ov\mu_i<0\right\}$.
\end{rem}

\subsection{Type $B_{n}$.}\label{tipoB}
\vskip10pt
We have $\Dp=\{e_i\pm e_j\mid 1 \leq i<j\leq n\}\cup \{e_i\mid 1\leq
i\leq n\},\r=\sum\limits_{i=1}^n\frac{2n-2i+1}{2}e_i,$ $\th=e_1+e_2$, hence 
$h^\vee=2n-
1$. Moreover
\begin{align*}
&P =\left\{\sum_{i=1}^n\frac{x_i}{2} e_i\mid x_i\text{ all even or
all odd}\right\},\qquad
Q =\sum_{i=1}^{n}\ganz e_i,
\\
&P^+=\left\{\sum_{i=1}^n\l_i e_i\in P\mid \l_1\geq\l_2\geq\dots\l_n\geq 
0\right\}.\end{align*}

Since $(\th,\th)=2$ we have $(\cdot,\cdot)=\frac{1}{2 
h^\vee}\langle\cdot,\cdot\rangle$. 

\begin{prop}\label{tipob} 
  For
$\l=\sum\limits_{i=1}^n\l_ie_i\in P^+$ we have
$$\l\in P_{alc}\iff \begin{array}{llr} \l_i\in\ganz\quad \text{for \ 
$i=1,\ldots,n$,}\\
2(\l_i+n-i)+1\not\equiv \pm 2(\l_j+n-j)+1
\mod\,2(2n-1)\\(i\ne j).\end{array}$$
\end{prop}
\begin{proof} By Corollary \ref{perkostant} we have that $\l\in 
Q=\sum\limits_{i=1}^n\ganz e_i$.
The second condition follows directly from \eqref{cond} and the observation that
$\langle \l+\r,e_i\rangle\notin\ganz$ for $ i=1,\ldots,n$.
\end{proof}
\vskip5pt Observe that 
\begin{align*}
\tfrac{1}{2}\Q&=\frac{1}{2}\left\{\tau\in\h^*_\R\mid 
(\tau,x)\in\ganz\,\forall\,x\in P^+\right\}\\
&=\frac{1}{2}\left\{\tau\in\h^*_\R\mid \langle\tau,x\rangle\in 2h^\vee 
\ganz\,\forall\,x\in P^+\right\}\\
&=h^\vee\left\{\tau\in\h^*_\R\mid \langle\tau,x\rangle\in \ganz\,\forall\,x\in 
P^+\right\}\\
&=h^\vee\{\tau=\sum_{i=1}^n\tau_ie_i\in\h^*_\R\mid
\tau_i\in\ganz,\,\sum_{i=1}^n\tau_i\text { even}\}.\end{align*}
Assume that $\l\in P_{alc}$  and set $\mu=\l+\r$, so that
$\mu=\sum\limits_{i=1}^n\frac{\mu_i}{2}e_i$ with $\mu_i\in 2\ganz + 1$ for 
$i=1,\ldots,n$. Denote by $\ov\mu_i$ the unique
element in $\{\pm 1,\pm 3,\ldots,\pm (2n-3)\}\cup\{2n-1\}$ such that 
$\mu\equiv\ov\mu_i\,\mod\,2(2n-1)$ and set
$\tilde\mu=\sum\limits_{i=1}^n\frac{\ov\mu_i}{2}e_i$. Consider $\mu-\tilde\mu$: 
if 
$\mu-\tilde\mu\in\frac{1}{2}\Q$
we set $\ov \mu=\tilde \mu$. Otherwise let $i^*$ be the unique index such 
that $\mu_{i^*}=2n-1$.
and set $\ov\mu=\tilde\mu-\frac{2n-1}{2}e_{i^*}$. This is equivalent to changing $2n-1$ 
into $-(2n-1)$ in the sequence
of  remainders. 
Then we obtain that $\mu-\ov\mu\in\frac{1}{2}\Q$. Now we observe that in any 
case  $\ov\mu\in W\cdot\r$,
say $\ov\mu=v(\r)$. Hence if we set $\tau=\mu-\ov\mu$, we obtain that 
$\mu=\l+\r=t_\tau v(\r)$ and 
$\chi_\l(a)=(-1)^{\ell(v)}$.
 
\begin{rem}\label{segnob} If $v\in W$, we define   $(\pm i)^{\s_v}=\pm2\left\langle
v(\r),e_{n-(i-1)/2}\right\rangle$ for $i=1,3,\dots,2n-1$. 
Since also in type $B$ the Weyl group
acts as signed permutations on 
$\{\pm e_1,\ldots,\pm e_n\}$ we have that the map $v\mapsto \s_v$ defines an 
homomorphism from $W$ to the set of signed
permutations on $\{\pm 1,\pm 3,\ldots, \pm (2n-1)\}$. Arguing as in type $C$ we find that
 $\chi_\l(a)=sign(\left|\s_v\right|)(-1)^{n_{\s_v}}$ where $|\s_v|$ is the
 permutation of $\left\{1,3,\dots,2n-1\right\}$ defined by setting
$i^{\left|\s_v\right|}=|\ov\mu_{n-(i-1)/2}|$ and $n_{\s_v}=\sharp\left\{i\mid
\ov\mu_i<0\right\}$.
\end{rem}
\vskip10pt
\subsection{Type $D_{n}$.}\label{tipoD}
\vskip10pt
We have $\Dp=\{e_i\pm e_j\mid  1\leq i<j\leq n\}, \r=\sum\limits_{i=1}^n(n-
i)e_i,$ 
$\th=e_1+e_2$, hence
$h^\vee=2n-2$. Moreover
\begin{align*}
&P =\left\{\sum_{i=1}^n\frac{\l_i}{2} e_i\mid \l_i\text{ all even or
all odd}\right\},\\
&Q =\left\{\sum_{i=1}^n\l_i e_i\mid \sum_{i=1}^n\l_i\,\text{even}\right\},
\\
&P^+=\left\{\sum_{i=1}^n\l_i e_i\in P\mid 
\l_1\geq\l_2\geq\dots\geq|\l_n|\right\}.\end{align*}
Since $(\th,\th)=2$ we have $(\cdot,\cdot)=\frac{1}{2 
h^\vee}\langle\cdot,\cdot\rangle$.
As in type $B_n$, Corollary \ref{perkostant} implies the following result.

\begin{prop}\label{tipod} 
  For
$\l=\sum\limits_{i=1}^n\l_ie_i\in P^+$ we have
$$\l\in P_{alc}\iff \begin{array}{llr} \l_i\in\ganz\quad \text{for \ 
$i=1,\ldots,n$,}\
\sum\limits_{i=1}^n\l_i\,\text{even,}\\
\l_i+n-i\not\equiv \pm(\l_j+n-j)
\mod\,(2n-2)\ (i\ne j).\end{array}$$
\end{prop}
\vskip5pt Observe that in this case
$
\frac{1}{2}\Q=h^\vee Q$.
Assume that $\l\in P_{alc}$  and set $\mu=\l+\r$, so that
$\mu=\sum\limits_{i=1}^n\mu_ie_i$ with $\mu_i\in \ganz$ for $i=1,\ldots,n$. 
Denote by 
$\ov\mu_i$ the unique
element in $\{\pm 1,\pm 2,\ldots,\pm (n-2)\}\cup\{0,n-1\}$ such that 
$\mu\equiv\ov\mu_i\,\mod\,(2n-2)$ and set
$\tilde\mu=\sum\limits_{i=1}^n\ov\mu_ie_i$. Consider $\mu-\tilde\mu$: if $\mu-
\tilde\mu\in\frac{1}{2}\Q$
we define $\ov \mu=\tilde \mu$. Otherwise let $i^*$ be the unique index such 
that $\mu_{i^*}=n-1$
and set $\ov\mu=\tilde\mu-2(n-1)e_{i^*}$. This is equivalent to changing $n-1$ 
into $-(n-1)$ in the sequence
of remainders. 
Then we obtain that $\mu-\ov\mu\in\frac{1}{2}\Q$. As in type $B_n$ we have 
 $\ov\mu=v(\r),\,v\in W$ and $\mu=\l+\r=t_\tau v(\r)$ with $\tau=\mu-\ov\mu$. As 
before, 
$\chi_\l(a)=(-1)^{\ell(v)}$.

\begin{rem}\label{segnod}This time the action of $W$ on $\r$ defines an homomorphism $v\mapsto |\s_v|$ onto the set of  permutations on $\left\{0,1,2,\dots,n-1\right\}$. The permutation $|\s_v|$ is defined by setting   $i^{|\s_v|}=|\left\langle v(\r),e_{n-i}\right\rangle|$. Since $|\s_{s_i}|$ is a simple transposition, it follows as before that $(-1)^{\ell(v)}=sign(|\s_v|)$, hence $\chi_\l(a)$ is the sign of the permutation of $\left\{0,1,\dots,n-1\right\}$ defined by setting $i\mapsto|\ov\mu_{n-i}|$.
\end{rem}
\vskip10pt

\vskip10pt
\subsection{Type $G_{2}$.}\label{tipoG}
\vskip10pt
It is amusing to work out our Euclidean algorithm for type $G_2$ also. Following 
\cite{bourbaki} we realize the root system of type $G_2$ in 
$$
V=\left\{(x_1,x_2,x_3)\in \R^3\mid x_1+x_2+x_3=0\right\}.
$$
 As above $\left\langle\cdot,\cdot \right\rangle$ is the standard inner product 
on $\R^3$ and $\left\{e_1,e_2,e_3\right\}$ is the canonical basis.
We have 
$$
\D=\{\pm( e_i- e_j)\mid  1\leq i,j\leq 3\}\cup\left\{\pm(2e_i-e_j-e_k)\mid 
\{i,j,k\}=\{1,2,3\}\right\},
$$
 $\Pi=\left\{e_1-e_2,-2e_1+e_2+e_3\right\},
$
 so that $\r=-e_1-2e_2+3e_3$, 
$\th=-e_1-e_2+2e_3$, hence
$h^\vee=4$. Moreover
\begin{align*}
&P=Q=V\cap(\sum_{i=1}^3\ganz e_i),
&P^+=\left\{\sum_{i=1}^3\l_i e_i\in P\mid 
0\ge\l_1\geq\l_2\right\}.\end{align*}
Since $(\th,\th)=6$ we have $(\cdot,\cdot)=\frac{1}{6 
h^\vee}\langle\cdot,\cdot\rangle$. Set $\varepsilon_i=-1$ for $i=1,2$ and 
$\varepsilon_3=1$.
Corollary \ref{perkostant} implies the following result.

\begin{prop}\label{tipog} 
  For
$\l=\sum\limits_{i=1}^3\l_ie_i\in P^+$ we have that $\l\in P_{alc}$ if and only 
if
\begin{align}\label{f1}
&\l_i+\varepsilon_ii\not\equiv \l_j+\varepsilon_jj
\!\!\!\!\!\!&&\mod\,(12)\ (i\ne j)\\
&2(\l_i+\varepsilon_ii)\not\equiv\l_j+\varepsilon_jj +\l_k+\varepsilon_kk\!\!\label{f2}
&&\mod\,(12)\,(\{i,j,k\}=\{1,2,3\})
\end{align}
\end{prop}

\vskip5pt 

An easy calculation shows  that in this case
\begin{equation}\label{Qg2}
\frac{1}{2}\Q=
4\left\{\sum\limits_{i=1}^3 x_i e_i\in Q\mid x_1\equiv x_2\equiv x_3\mod 
(3)\right\}.
\end{equation}
Assume that $\l\in P_{alc}$  and set $\mu=\l+\r$, so that
$\mu=\sum\limits_{i=1}^3\mu_ie_i$ with $\mu_i=\l_i+\varepsilon_ii\in \ganz$ and 
$\mu_1+\mu_2+\mu_3=0$.
 Denote by 
$[\mu_i]_n=\mu_i+n\ganz\in\ganz/n\ganz$. By the chinese remainder theorem the 
map
$[\mu_i]_{12}\mapsto([\mu_i]_3,[\mu_i]_4)$ is an isomorphism.

Since $\sum\limits_{i=1}^3\mu_i=0$, we have obviously that 
$\sum\limits_{i=1}^3[\mu_i]_n=0$.
Relation \eqref{f1} implies that $([\mu_i]_3,[\mu_i]_4)\ne([\mu_j]_3,[\mu_j]_4)$ 
if $i\ne j$. Moreover we have the following further conditions:	
\begin{equation}\label{condizionig}
	\begin{array}{ll}
		\text{$\left[\mu_j\right]_4$ cannot be all equal,}\\
	\left[\mu_j\right]_4\ne0&j=1,2,3,\\
	\left[\mu_i\right]_4+\left[\mu_j\right]_4\ne 0\,	
	&\text{if $i\ne j$.}
	\end{array}
\end{equation}
Let us check  the first condition: if $[\mu_1]_4=[\mu_2]_4=[\mu_3]_4=x$ then
$$
-
2([\mu_1]_3,x)+([\mu_2]_3,x)+([\mu_3]_3,x)=([\mu_1]_3+[\mu_2]_3+[\mu_3]_3,0)=(0,
0)
$$
and this contradicts  \eqref{f2} . For the second condition suppose
$[\mu_i]_4=0$. Let $j,k$ be such that $\left\{i,j,k\right\}=\left\{1,2,3\right\}$. Since $[\mu_i]_4+[\mu_j]_4+[\mu_k]_4=0$ we have that $-2[\mu_i]_4+[\mu_j]_4+[\mu_k]_4=-3[\mu_i]_4=0$ hence
$$
-
2([\mu_i]_3,[\mu_i]_4)+([\mu_j]_3,[\mu_j]_4)+([\mu_k]_3,[\mu_k]_4)=([\mu_i]_3+[\mu_j]_3+[\mu_k]_3,0)=(0,
0).
$$
The third condition is obtained in the same way.

Set $S=\left\{([\mu_i]_3,[\mu_i]_4)\mid i=1,2,3\right\}$. The conditions in 
\eqref{condizionig} 
imply that there are two possibilities for $S$: either 
$S=\left\{(a,[1]_4),(b,[1]_4),(c,[2]_4)\right\}$
or $S=\left\{(a,[3]_4),(b,[3]_4),(c,[2]_4)\right\}$. Relation \eqref{f1} forces 
$a\ne b$, so that
$a-b=\pm[1]_3$.  Define the ordered sets
\begin{align*}
&S_1=((a,[1]_4),(b,[1]_4),(c,[2]_4)),\\
&S_2=((a,[3]_4),(b,[3]_4),(c,[2]_4)).\end{align*}
The algorithm works as follows. Let $i^*,j^*,k^*$ be such that 
$([\mu_{i^*}]_{12},[\mu_{j^*}]_{12},$ $[\mu_{k^*}]_{12})=S_x,\,x=1,2$, and write
$\mu_y=4\tilde q_y+\tilde r_y,\,y\in\{i^*,j^*,k^*\}$, where the sequence of 
remainders
$\tilde r_y$ is $(1,1,2)$ if $x=1$ and $(3,3,2)$ if $x=2$; this of course 
determines the $\tilde q_y$.
Now change the sequence of quotiens $\tilde q_y$ into a new sequence $q_y$ in  
such a  way to obtain 
the following new remainders $r_y$
$$\begin{array}{lllll}
&&&&(r_{i^*},r_{j^*},r_{k^*})\\
x=1&\quad&a-b=[1]_3&\quad&(1,-3,2)\\
x=1&\quad&a-b=-[1]_3&\quad&(-3,1,2)\\
x=2&\quad&a-b=[1]_3&\quad&(3,-1,-2)\\
x=2&\quad&a-b=-[1]_3&\quad&(-1,3,-2)
\end{array}$$
This choice implies 
$q_{i^*}\equiv q_{j^*}\equiv q_{k^*}\mod (3)$. For instance assume 
$x=1,a-b=[1]_3$.  Since $a=[q_{i^*}+1]_3$, $b=[q_{j^*}]_3$, and $c=[q_{k^*}+2]_3$,  we 
have that
$0=a-b-[1]_3=[q_{i^*}-q_{j^*}]_3$ and, since $ \sum\limits_{i=1}^3q_i=0$ we also 
obtain that
$[q_{i^*}-q_{k^*}]_3=0$. The other cases are checked similarly.
\par
In all cases we have that, if we set $\tau=\sum_i q_ie_i$ then $\tau\in\half \Q$. Moreover 
$\left\{r_1,r_2,r_3\right\}=\pm\left\{1,2,-3\right\}$. 
We now observe  that  $\sum_i r_i e_i$ is in $W\cdot\r$. This is an immediate 
consequence of the general fact that, if $\l\in P$ and $(\l,\l)=(\r,\r)$, then 
$\l=w\r$ for some $w\in W$. (A less attractive proof is obtained by simply 
listing all twelve elements of $W\cdot\r$).
Thus $\mu=t_\tau v(\rho)$, where $v$ is the unique element of  $W$ such that 
$v(\rho)=\sum_i 
r_ie_i$.\par
A more explicit description of $v$ and the determination of $\chi_\l(a)$ will be 
performed 
at the end of Section 4.


\section
{Affine Weyl groups as permutations of $\ganz$.}\label{permutazionitutti}

In this section we will show how one can construct realizations of the
classical affine Weyl groups as permutations of $\ganz$ from the knowledge of 
the orbit $\Wa_q\cdot\l$, for an appropriate choice of $\l$ and $q$.
Our treatment takes into account all the representations of classical affine 
Weyl 
groups known in literature.
We  obtain analogous results also for $\tilde{G_2}$.

\vspace{5pt}

We shall use several times the following obvious facts.

\begin{vuoto}
\label{ext}
Let $p\in\nat^+$ and assume that:
\item (1) $A=\{a_1, \dots, a_p\}$ is a set of representatives of $\ganz/p\ganz$; 
\item (2) $f:A\to \ganz, \ a_i\mapsto a_i^f$ is a map such that 
$\{a_1^f, \dots, a_p^f\}$ is still a set of representatives  of $\ganz/p\ganz$. 
\par
Then $\tilde f:\ganz \to \ganz, \ a_i+kp\mapsto a_i^f+kp$ for all $k\in 
\ganz$, is a permutation of $\ganz$ which extends $f$.
\end{vuoto}

\begin{vuoto}
\label{ovvio}
Let $q\in\R^+$ and assume that $\l\in\h^*_{\mathbb R}$ is $q$-regular. 
Then $w\mapsto w(\l)$ is a bijection from $\Wa_q$ to the orbit $\Wa_q\cdot 
\l$ of $\l$ under $\Wa_q$. 
\end{vuoto}

\vspace{5pt}
\par
{\bf \large Types $\tilde A_{n-1}$,  $\tilde C_n$,  
$\tilde B_n$, and  $\tilde D_n$.} 

\vspace{5pt}

We shall use the following notation: 
for $a,b\in \ganz$ with $a<b$, $c\in\ganz$ with $c>0$, $A\subseteq \ganz$ we 
set 
$$
[a,b]=\{z\in \ganz\mid a\leq z\leq b\},\quad [c]=[1,c]; \quad \pm A=A\cup-A.
$$
For any set $N$, we denote by $S(N)$ the group of permutations of $N$.

\vspace{5pt}

We realize the classical root systems as in \cite{bourbaki}, except that we
reverse the order of the canonical basis of $\real^n$.  
Thus if $\{e_i\mid i\in [n]\}$ is the canonical basis of $\real^n$,  
the simple roots and the highest root are:

\vspace{5pt}

\noindent
for $A_{n-1}$:\quad  $\a_i=e_{i+1}-e_i$ for $i=1,\ldots,n-1$;
\quad $\th=e_{n}-e_1$; 
\par\noindent
for $C_n$: \quad $\a_1=2e_1,\,\a_i=e_{i}-e_{i-1}$ for $i=2,\ldots,n$;
\quad $\th=2e_n$; 
\par\noindent
for $B_n$:\quad $\a_1=e_1$, $\a_i=e_i-e_{i-1}$ for $i=2, \dots, n$;
\quad $\th=e_{n-1}+e_n$;   
\par\noindent
for $D_n$:\quad $\a_1=e_1+e_2$, $\a_i=e_i-e_{i-1}$ for $i=2, \dots, n$;
\quad $\th=e_{n-1}+e_n$. 

\vspace{5pt}

If $\D$ is of type $A_{n-1}$, then $\D$ is a subset of $V= 
\{\sum\limits_{i=1}^{n}x_ie_i 
\mid \sum\limits_{i=1}^{n}x_i=0\}$.
We extend the faithful action of $W$ on $V$ to $\R^{n}$ by fixing 
pointwise $V^\perp$. We also naturally extend the translation action 
of $\Wa_q$ to $\real^n$.

\vspace{5pt}

Set 
$$
\l=\sum\limits_{i\in[n]} ie_i.
$$ 
Observe that 
$$
\l=
\begin{cases}
n\l_0+\ov\r\quad&\text{ in type $A_{n-1}$,}\\
\r\quad&\text{ in type $C_n$,}\\
\r+\o_1\quad&\text{ in type $B_n$,}\\
\r+2\o_1\quad&\text{ in type $D_n$.}
\end{cases}
$$

We set 
$\Qst=\sum\limits_{\a\in\Pi}\ganz\,\frac{2\a}{\langle\a,\a\rangle}$, thus

$$
\Qst=\frac{1}{c}\Q
$$ 
with $c=\langle\th,\th\rangle h^\vee$. 
The element $\l$ is $\frac {p}{c}$--regular where

$$
p=
\begin{cases}
n \quad&\text{ in type $A_{n-1}$,}\\
2n+1 \quad&\text{ in types $B_n$, $C_n$, and $D_n$.}
\end{cases}
$$
In particular, by Fact \ref{ovvio}, $w\mapsto w(\l)$ is a bijection from 
$\Wa_\psuc$ to $\Wa_\psuc\cdot\l$. 
We notice that 

$$
\Wa_\psuc= p\Qst\rtimes W,
$$
where we identify $p\Qst$ with the group of translations of $\real^n$
by elements of $p\Qst$. 
We also observe that for types $A_n$ and $C_n$ we have $\frac{p}{c}=\half$.
\par\vskip5pt
We set 
$$
I=
\begin{cases}
[n]   \quad&\text{ in type $A_{n-1}$ },\\
\left[-n,n\right] \quad&\text{ in types $B_n$, $C_n$, $D_n$}.\\
\end{cases}
$$
Thus $I$ is a set of representatives of $\ganz/p\ganz$. 
For types $B_n$, $C_n$, and $D_n$, we set 

$$
e_0=0, \quad e_{-i}=-e_{i}
$$ 
for all $i\in [n]$. Thus $e_i$ is defined for all $i\in I$. It is well-known 
that the finite Weyl group $W$ permutes $\{e_i \mid i\in I\}$. 
\par 

For all $w\in \Wa_\psuc$, and $i\in I$, we set 

\begin{equation}
\label{star}
\iw=\lan w(\l),e_i\ran.
\end{equation}
Then, by Fact \ref{ovvio}, $\ww$ determines $w$. Since $\langle\cdot,\cdot\rangle$ 
is $W$-invariant and $W$ permutes the $e_i$, for $w\in W$ we have that  

$$
e_{\iw}=w^{-1}e_i.
$$
This makes clear that $w\mapsto \ww$ is an injective homomorphism of the finite 
Weyl group $W$ into $S(I)$. In fact, this is the usual 
permutation representation of $W$. For $A_{n-1}$, $\{\ww\mid w\in W\}$ is the 
whole 
symmetric group $S_n$; for both $C_n$ and $B_n$, $\{\ww\mid w\in W\}$ is the 
group of all permutations of $[-n,n]$ such that $(-i)^{\ww}=-\iw$; for $D_n$, 
$\{\ww\mid w\in W\}$ is the group of all permutations of $[-n,n]$ such that 
$(-i)^{\ww}=-\iw$ and $|\{i\in [n]\mid \iw<0\}|$ is even. 
\par

We recall that for type $A_{n-1}$ the lattice $\Qst$ is the subgroup of 
$\sum\limits_{i\in[n]}\ganz e_i$ with zero coordinate sum. 
For type $C_n$, $\Qst=\sum_{i\in[n]}\ganz e_i$, while for both $B_n$ and 
$D_n$, $Q^\vee$ is the subgroup of $\sum\limits_{i\in[n]}\ganz e_i$ of all 
elements with even coordinate sum. In particular, since $\Wa_{\frac 
{p}{c}}=p\Qst\rtimes W$, we obtain in any case that for all $w\in 
\Wa_{\frac {p}{c}}$ and $i\in I$

$$
\iw\in\ganz \quad\text{and}\quad \{\iw\mod p\mid i\in I\}=\{i \mod p\mid 
i\in I\}.
$$
Thus, since $I$ is a set of representatives of $\ganz/p\ganz$, the map 
$\ww$ satisfies conditions (1) and (2) of Fact \ref{ext}. It follows that $\ww$ 
extends to a bijection of $\ganz$ onto itself, which we still denote by 
$\ww$, defined by 
\begin{equation}
\label{periodo}
(i+kp)^\ww=\iw+kp
\end{equation}
for all $i\in I$. We notice that in types $C_n$, $B_n$, and $D_n$, 
since $0^{w_*}=0$, we have that
$
z^{w_*}=z$ for all
$
z\in p\ganz=(2n+1)\ganz$ and $w\in \Wa_\psuc$.
\par

We shall verify that $w\mapsto \ww$ is an injective homomorphism of 
$\Wa_\psuc$ into the group of all permutations of $\ganz$. 
It is obvious that $\ww$ is uniquely determined by $I^\ww$, and hence by 
$w(\l)$,  
so injectivity follows immediately from Fact \ref{ovvio}.
Assume $\wa, \ua\in \Wa$, $\wa=t_\eta w$, $\ua=t_\tau u$, with $w, u\in W$ and 
$\tau, \eta\in p\Qst$. Then for $i\in I$
$$
i^{\wa_*}=\lan \wa(\l),e_i\ran =\lan \eta, e_i\ran+\lan w(\l), e_i\ran
=\lan \eta, e_i\ran +i^\ww,
$$
and since $\lan \eta, e_i \ran \in p\ganz$ and $ e_{\iw}=w^{-1}(e_i)$, we obtain 
$$
(i^{\wa_*})^{\ua_*}=\lan \eta, e_i \ran +i^{\ww\ua_*}=\lan \eta, e_i \ran +
\lan \tau,e_{\iw}\ran +i^{\ww 
u_*}=\lan \eta+w(\tau), e_i \ran +i^{\ww u_*}.
$$

On the other hand, $\wa\ua=t_{\eta+w(\tau)}wu$, hence  
$$
i^{(\wa\ua)_*}=\lan \eta+w(\tau), e_i \ran +i^{(wu)_*},
$$
and since $i^{(wu)_*}=i^{\ww u_*}$, we obtain that 
$i^{(\wa\ua)_*}=i^{\wa_*\ua_*}$.
\par 

\begin{rem}Suppose that we are given a homomorphism $w\mapsto w'$ from $\Wa_\psuc$ to $S(\ganz)$ such that \eqref{star} holds. Then,
for all $w\in \Wa_\psuc$,  
$w'_{|I}=w_{*|I}$. If  $w\in \Wa_\psuc$ and $u=t_\eta$, $\eta\in p\Qst$, 
then, for $i\in I$,
\begin{align*}
i^{(uw)'}&=\lan \eta+w(\l), e_i\ran= 
\lan \eta, e_i\ran +i^{w'}\\
&=i^{u'w'}=
\lan \eta+\l, e_i\ran^{w'}=(\lan \eta, e_i\ran+i)^{w'}.
\end{align*}
From the explicit description of $\Qst$, it is clear that for 
all $i\in I$ and $k\in \ganz$ there exists $\eta\in \Qst$ such
that $\lan \eta, e_i\ran=k$. It follows that relation \eqref{periodo} 
holds with  $w'$ in place of $w_*$, and therefore $w'=w_*$.
Thus the $w_*$ are the only permutations of $\ganz$ such that \eqref{star}
 holds and $w\mapsto w_*$ is a homomorphism of $\Wa_\psuc$
into $S(\ganz)$. 
\end{rem}

Combining the previous discussion with the results of Section \ref{orbita} we 
obtain Lusztig's description of the affine group of type $\tilde A_{n-1}$ 
\cite[\S\ 3.6]{L}. Recall that, in this case, $p=n=h^\vee$, 
and $\psuc=\half$.

\begin{theorem}
\label{A} 
If $\D$ is of type $A_{n-1}$, the map $w\mapsto \ww$ is a permutation 
representation of $\Wa_\half$ in $S(\ganz)$. Its image $\{w_*\mid 
w\in\Wa_\half\}$ is the group of all $f\in S(\ganz)$ such that 
\item 
(1) $(z+n)^f=z^f+n$ for all $z\in \ganz$; 
\item 
(2) $\sum\limits_{i=1}^{n} i ^f=\sum\limits_{i=1}^{n}i$.
\end{theorem}

\begin{proof} 
The first statement has already been proved.
It is clear from definitions that $(z+n)^{w_*}=z^{w_*}+n$ for 
all $w\in \Wa_\half$. 
It is also clear that condition (2) holds for all $v\in W$. 
If $w\in \Wa_\psuc$, $w=t_\eta v$, $\eta\in n\Qst$, $v\in W$, 
then  
$$
\sum_{i=1}^{n}i^{w_*}
=\sum_{i=1}^{n}\langle\eta,e_i\rangle+\sum_{i=1}^{n} i ^{v_*}.
$$
But, it is obvious, by the explicit description of $\Qst$, 
that $\sum\limits_{i=1}^{n}\langle\eta,e_i\rangle=0$, hence (2) holds for $w$. 
\par
It remains to prove that if $f\in S(\ganz)$ satisfies (1) and (2), 
then there exists $w\in \Wa_{\half}$ such that $f=w_*$. Let
$f$ be such that (1) and (2) hold and set 
$a_i=i^f$, $i=1,\dots,n$. Then
$a_i\ne a_j\mod  n$ if $i\ne j$ (otherwise $f$ is not a bijection). 
It follows from Proposition 
\ref{fund} that
$$
\sum_{i=1}^{n}\left(a_i-\frac{1}{n}\sum_{j=1}^{n}a_j\right)e_{i}=w(\r)
$$ 
for some $w\in\Wa_\half$. Observe that $\frac{1}{n}\sum\limits_{j=1}^{n}a_j=
\frac{n+1}{2}$, hence
$\sum\limits_{i=1}^{n}a_ie_i=\frac{n+1}{2}\l_0+w(\r)=w(\l)$. This implies that
$\langle w(\l),e_i\rangle=a_i$, hence $f=\ww$.
\end{proof}
\begin{rem}
The  affine reflection $s_0$ is equal to 
$t_{\frac {\th^\vee}{2}}s_\th$. Since $\th=e_n-e_1$ and 
$\frac{\th^\vee}{2}=n(e_n-e_1)$, we obtain that
\begin{align*}
j^{s_{0*}}&=\lan t_{n(e_n-e_1)} s_\th(\l),e_j\ran \\
&=\lan
n(e_n-e_1)+ne_1+\sum_{i=2}^{n-1}ie_i+e_n ,e_j\ran=
\begin{cases} 
0\ \text{for}\ j=1,\\
j\ \text{for}\ 2\leq j\leq n-1,\\
n+1 \ \text{for}\ j=n.\\
\end{cases}
\end{align*}

Clearly, for $i\in[n-1]$, $s_{i*}$ acts on $[n]$ as the
transposition $(i, i+1)$. 
\end{rem}

\begin{rem} We may apply formula \eqref{alcove} with $\mu=\l$. 
Since positive roots in $A_{n-1}$ are of the form $\a_{ij}=e_j-e_i,\,i<j$, we
deduce, using \eqref{killinga}, the following relation

$$
k(w,\a_{ij})=
\left\lfloor\frac{(w(\l),\a_{ij})}{\half}\right\rfloor=
\left\lfloor\frac{\langle w(\l), e_j-e_i\rangle}{h^\vee}\right\rfloor=
\left\lfloor\frac{j^{w_*}-i^{w_*}}{n}\right\rfloor.
$$
This is one statement of Theorem 4.1 from \cite{Shi} (taking into account the 
different notational conventions). We also have, by \eqref{lunghezza} 
$$
\ell_{\half}(w)=\sum_{1\leq i<j\leq 
n}|\left\lfloor\frac{j^{w_*}-i^{w_*}}{n}\right\rfloor|,
$$ 
a formula which appears, with different derivations, in 
\cite{BB},\cite{E},\cite{Papi},\cite{Shi}.
\end{rem}

\vspace{5pt}

\begin{theorem}
\label{C}
If $\D$ is of type $C_n$, then $w\mapsto \ww$ is an injective homomorphism 
of $\Wa_\half$ into $S(\ganz)$. 
Its image $\{\ww\mid w\in \Wa_\half\}$ is the subgroup of all 
permutations $f$ of $\ganz$ such that
\item
(1) $(-z)^f=-z^f$ for all $z\in \ganz$; 
\item 
(2) $(z+k(2n+1))^f=z^f+k(2n+1)$ for all $z, k\in \ganz.$
\end{theorem}

\begin{proof}
Recall that in this case $p=2n+1$ and $\psuc=\half$.
It follows directly from  definitions that, for all $w\in \Wa_\psuc$, $\ww$ 
satisfies conditions (1) and (2).
It remains to prove that all permutations of $\ganz$ which satisfy conditions 
(1) and (2) lie in $\{\ww\mid w\in \Wa_\psuc\}$.
\par

The anti-symmetry condition (1) implies in particular that $0^f=0$, hence any 
odd $f\in S(\ganz)$ satisfies (2) if and only if it permutes the non zero 
cosets in $\ganz/p\ganz$. This means that $\{0, \pm 1^f, \dots ,\pm n^f\}$
is a set of representative of $\ganz/p\ganz$ or, equivalently, that 

\begin{equation}
\label{effe}
i^f\not\equiv 0 , \ i^f\pm j^f\not\equiv 0\mod p,\quad\text 
{for}\quad  1\leq i < j \leq n
\end{equation}
(notice that $p$ being odd,  $i^f\not\equiv 0\mod p$ if and only if 
$2i^f\not\equiv 0 \mod p$).
\par

Now we recall that $P=\Qst$ and $\l=\rho$, so that, by Lemma \ref{taut}, 
$\l+\Qst\cap C_\psuc=\{\l\}$. Since $\Wa_\psuc= p\Qst\rtimes W$, it is clear 
that $\Wa_\psuc$ acts on $\l+P$. By Proposition \ref{fund} we obtain that 
$\Wa_{\psuc}\cdot\l$ is the set of all $\mu\in \l+ P$ such that $(\mu, 
\a)\not\in\psuc \ganz$ or, equivalently, $\lan \mu, \a\ran \not\in p\Qst$ for 
each root  $\a$. By the explicit description of the root system, this 
means that, if $\mu=\sum\limits_{i=1}^n\mu_i e_i$, then
$$
2\mu_i,\ \mu_i\pm \mu_j\not\in p\ganz\quad\text{for}\quad 1\leq i< j\leq n.
$$
Comparing the above conditions with \eqref{effe}, we deduce that for each 
$f\in S(\ganz)$ such that (1) and (2) hold, there exists $w\in \Wa_\psuc$ 
such that $\sum\limits_{i=1}^ni^f e_i= w(\l)$, and therefore such that $f=\ww$.
\end{proof}

\begin{rem}
In our setting, the  affine reflection $s_0$ is equal to 
$t_{\frac {2n+1}{c}\th^{\vee}}s_\th$. Since $\th=2e_n$ and 
$\frac {1}{c}\th^{\vee}=\half \th=e_n$, we obtain that
$$
j^{s_{0*}}=\lan t_{(2n+1)e_n} s_\th(\l),e_j\ran =
\lan e_n+\l,e_j\ran =
\begin{cases} 
j\ &\text{for}\ 1\leq j< n,\\
n+1 \ &\text{for}\ 1\leq j< n.\\
\end{cases}
$$

Clearly, for $i\in[n-1]$, $s_{i*}$ acts on $[n]$ as the
transposition $(i, i-1)$, while $s_{n*}$ acts on $[-n,n]$
as the transposition $(-n, n)$. 
\end{rem}

\begin{rem}
The representation of the Weyl group of type $\tilde C_n$ 
as a subgroup of $S(\ganz)$ obtained in Therem \ref{C} coincides with 
the one presented by Bedard \cite{Bedard}. 
A different representation appears in literature (see \cite{Shi}, \cite{Papi}). 
We can as well get this representation in our framework. Indeed, we note that, 
with the notation of Lemma \ref{taut}, there are two possible values of $q$ 
verifying equation \eqref{eqtaut}: $2n+1$ and $2n+2$. Hence we can  
define an injective homomorphism $w\mapsto w_{**}$  
of $\Wa_{\frac {2n+2}{c}}$ into $S(\ganz)$ setting 
\begin{align*}
&i^{w_{**}}=\lan w(\l), e_i\ran\ 
\text{for}\ i\in [-n,n],\ \pm(n+1)^{w_{**}}=\pm(n+1),\\
&(i+k(2n+2))^{w_{**}}=i^{w_{**}}+k(2n+2).
\end{align*}
Then $s_{i**}$ and $s_{i*}$ have the same action on $[-n, n]$, for $i\in [n]$.
The action of $s_{0**}$ is defined by 
$j^{s_{0**}}=j$ for $1\leq j< n$,  $n^ {s_{0*}}=n+2$, and by the  
condition of compatibility with translation by  $2n+2$.
\end{rem}

\vspace{5pt}

\begin{theorem}\label{B}
If $\D$ is of type $B_n$ then $w\mapsto \ww$ is an injective homomorphism 
of $\Wa_\psuc$ into $S(\ganz)$. 
Its image $\{\ww\mid w\in \Wa_\psuc\}$ is the subgroup of all 
permutations $f$ of $\ganz$ such that
\item
(1) $(-z)^f=-z^f$ for all $z\in \ganz$; 
\item 
(2) $(z+k(2n+1))^f=z^f+k(2n+1)$ for all $z, k\in \ganz$;
\item 
(3) $\sum\limits_{i=1}^ni^f\equiv {\binom {n+1}{2}}\mod 2$.
\end{theorem}

\begin{proof}
It remains to prove that $w_*$ satisfies (3) for all $w\in \Wa_\psuc$ 
and that each $f\in S(\ganz)$ such that (1), (2), and (3) hold is equal to some
$w_*$, $w\in \Wa_\psuc$.
\par

If $w\in W$, then $\{1^\ww, \dots, n^\ww\}$ differs from $[n]$ at most 
in the sign of elements, hence it is clear that 

$$
{\binom{n+1}{2}}=\sum_{i=1}^ni \equiv\sum_{i=1}^ni^\ww 
\mod 2.
$$
Since $\Qst$ is the the set of all elements in $\sum\limits_{i=1}^n\ganz e_i$ with 
even coordinate sum and $\Wa_\psuc=p\Qst\rtimes W$, it is clear that 
$w_*$ satisfies (3) for all $w \in \Wa_\psuc$.
\par

The above argument also shows that $\l+\Qst$ is $\Wa_\psuc$-stable.
Moreover, it is easily seen that  $\l+\Qst\cap C_\psuc=\{\l\}$. 
Thus we may apply Proposition \ref{fund}, with $L=\Qst$, so as to obtain 
that $\Wa_\psuc\cdot \l$ is the set of all $\mu\in \l+\Qst$ such that 
$\lan \mu, \a\ran \not\in p \ganz$, for each root $\a$.
From the explicit description of $\Qst$ and of the root system, we obtain
that,  if $\mu =\sum\limits_{i=1}^n \mu_i e_i\in \sum\limits_{i=1}^n\ganz
e_i$, then 
$\mu\in \Wa_\psuc\cdot \l$ if and only if 

$$
\sum_{i=1}^n\mu_i\equiv {\binom{n+1}{2}}\mod 2,\quad\text{and} 
\quad \mu_i,\ \mu_i\pm \mu_j\not\in p\ganz\quad\text{for}\quad 1\leq i< j\leq 
n.
$$

Now it is clear that the same argument used in the proof of Theorem \ref{C}
shows that if $f\in S(\ganz)$ satisfies condition (1), then condition (2) is 
equivalent to \eqref{effe}. 
We easily conclude that each $f\in S(\ganz)$ such that (1), (2), and (3)
hold is equal to $w_*$ for some $w\in\Wa_\psuc$.
\end{proof}
\begin{rem}\label{alt3}
Condition (3) in Theorem \ref{B} can be replaced by the following 
one:
\item
(3$'$) $\sum\limits_{i=1}^n \left(i^f-\ov{i^f}\right)\in 2(2n+1)\ganz$.
\par
or, equivalently,
\item
(3$''$) $|\{i\leq n\mid i^f>n\}|$ is even.
\par
In fact, if we set $i^f=k_i(2n+1)+\overline{i^f}$, then we 
have $\sum\limits_{i=1}^n \ov {i^f}\equiv \binom{n+1}{2}$ hence 
$\sum\limits_{i=1}^n {i^f}\equiv \binom{n+1}{2}$ if and only if 
$\sum\limits_{i=1}^n k_i$ is even, which is equivalent to 
condition (3$'$).
Moreover, $\{j\leq n\mid j^f>n\}=\sum\limits_{i=1}^n |k_i|\equiv \sum\limits_{i=1}^n k_i$,
and since $k_i(2n+1)=i^f-\overline{i^f}$, we obtain that
(3$''$) is equivalent to (3$'$) and hence to (3).
\end{rem}

\vspace{5pt}

We finally deal with type $D_n$. In this case, we identify 
$\Wa_\psuc$ with a subgroup of its $\tilde B_n$-analog. Namely,
if $W_{B_n}$ is the finite Weyl group for type $B_n$,  
we may identify the finite Weyl group of $D_n$ 
with the subgroup  of $W_{B_n}$

$$
W'=\{w\in W_{B_n}\mid i^{w_*}< 0 \ \text{for an even number of } i\in [n]\}
$$
and we set 

$$
\Wa_\psuc=p\Qst\rtimes W'.
$$

For $j\in \ganz$ we denote by $\ov j$ its residue modulo $p$. It 
is clear that  if $w\in\Wa_\psuc$, $w=t_\eta v$, with 
$\eta\in p\Qst$ and $v\in W'$, then $i^{v_*}=\ov{i^{w_*}}$ for all $i\in[n]$, 
and $\eta=\sum\limits_{i=1}^n(i^{\wa_*}-\iw)e_i$, hence from Theorem \ref{B} we 
directly obtain the following result.

\begin{theorem}
\label{D}
If $\D$ is of type $D_n$ then $w\mapsto \ww$ is an injective homomorphism 
of $\Wa_\psuc$ into $S(\ganz)$. 
Its image $\{\ww\mid w\in \Wa_\psuc\}$ is the subgroup of all 
permutations $f$ of $\ganz$ such that
\item
(1) $(-z)^f=-z^f$ for all $z\in \ganz$; 
\item 
(2) $(z+k(2n+1))^f=z^f+k(2n+1)$ for all $z, k\in \ganz$
\item 
(3) $\sum\limits_{i=1}^ni^f\equiv {\binom{n+1}{2}}\mod 2$, and 
 $|\{i\in [n]\mid \ov{i^f}<0\}|$ is even. 
\end{theorem}

\begin{rem}
For both types $B_n$ and $D_n$ we find that $s_0=t_{(2n+1)\th}s_\th$ and 
hence 
$$
s_0(\l)=(2n+1)\th+\l-\lan \l, \th\ran \th=\l+2\th=
\sum\limits_{i=1}^{n-2} i e_i+(n+1) e_{n-1}+(n+2)e_n.
$$
It follows that 
$$
i^{s_{0*}}=i \text{ \ for \ } i\in [n-2],\quad (n-1)^{{s_0}_*}=n+1,\quad 
n^{s_{0_*}}=n+2.
$$
Since $n+1=-n+(2n+1)$, and $n+2=-(n-1)+(2n+1)$, we have that $(n+1)^{{s_0}_*}=n-1$, and 
$(n+2)^{{s_0}_*}=n$. Thus ${s_0}_*$ acts on $\{-n+2, \dots, n+2\}$ as the 
product of transpositions $(n-1,n+1)(n,n+2)$.
For $i\in [n]$, the action of ${s_i}_*$ on $[-n,n]$ is the usual one, 
hence, for $2\leq i\leq n$, ${s_i}_*$ is the product of transpositions 
$(i-1,i)(-(i-1),-i)$; 
${s_1}_*$ is the transposition $(1,-1)$ for $B_n$, while 
is the product of transpositions $(1,-2)(2,-1)$ for $D_n$.
\par  
\end{rem}

\vspace{5pt}

{\bf \large Type $G_2$.}

\vspace{5pt}
In this case we shall define an injective homorphism of $\Wa$ ($=\Wa_1$) 
into $S(\ganz)$. We omit everywhere the subscript $1$, so $T$ is the
subgroup of translations of $\Wa$ and $C$ is the fundamental alcove. 
The rest of notation is the same of Subsection \ref{tipoG}.
The map $w\mapsto w_*$, $\Wa\to S(\ganz)$, we are going to define is 
determined by $w(\r)$. Injectivity will be an immediate consequence of
the fact that $\r\in C$. 
\par
We set $e_{-i}=-e_i$ for $i\in [3]$,  $\varepsilon_i=-1$ for $i=\pm 1,\pm2$, 
$\varepsilon_3=\varepsilon_{-3}=1$. Then we define, for all $w\in \Wa$, 
$$
0^{\ww}=0, \quad \iw=\varepsilon_i\lan w(\rho),e_i\ran\text{ \rm \ for \ } i\in 
\pm[3].
$$
If $v\in W$, and $i\in \pm[3]$, then there exist unique $j\in \pm[3]$ and 
$v_i\in V^\perp$ such that $v(e_i)=e_j+v_i$. Then for $w=v^{-1}$ we have  
$\lan w(\rho),e_i\ran= \lan\rho,v(e_i)\ran=\varepsilon_j j$, hence 
$$
w^{-1}(e_i)=\varepsilon_i\varepsilon_{\iw}e_{\iw}+v_i,
$$
with $v_i\in V^\perp$. It follows directly that for all $w,w'\in W$, 
$(ww')_*=w_*w'_*$, hence $w\mapsto w^*$ is an injective homomorphism 
of $W$ into the set of all permutations of $[-3,3]$.

 It is easily seen that 
the image $W_*$ of $W$ under this homomorphism is the set (group) of all 
permutations $f$ of $[-3,3]$ such that $(-i)^f=-i^f$ and $\sum\limits_{i\in [3]}
\varepsilon_i i^f=0$. Notice that this last condition is equivalent to 
$\{-1^f,-2^f,3^f\}$ being equal to either $ \{-1,-2,3\}$ or $\{1,2,-3\}$. By 
restricting maps
to $\pm[3]$ we obtain that the map $w\mapsto w_*$ defines an isomorphism between 
$W$ and the group 
of functions $f:\pm[3]\to\pm[3]$ such that $(-i)^f=-i^f$ and $\{-1^f,-2^f,3^f\}= 
\pm\{1,2,-3\}$.

We recall that 
$$
\Q=8\left\{\sum_{i=1}^3 x_i e_i\in Q\mid x_1\equiv x_2\equiv x_3\mod 3\right\}, 
$$
in particular, for each $t\in T$ and $i\in \pm[3]$, $i^{t_*}\equiv i\mod 8$. 
For all $w\in \Wa$, we define $4^\ww=4$. Then it is clear that 
$\ww$ maps the set of representatives $[-3,4]$ of $\ganz/8\ganz$ into some
set of representatives of $\ganz/8\ganz$, hence Fact \ref{ext} applies and 
$\ww$ can be extended to a bijection $\ww$ of $\ganz$ onto itself by 
setting $(i+8k)^\ww=i^\ww+8k$ for all $k\in \ganz$. Notice that $\ww$ fixes 
pointwise $4\ganz$.\par
We next verify that $w\mapsto \ww$ is an injective 
homomorphism of the whole $\Wa$ into the group of all permutations of $\ganz$. 
It is obvious that $\ww$ is determined by $[-3,3]^\ww$, hence by $w(\r)$, so, 
as remarked above, injectivity is immediate. Assume $\wa\in \Wa$, $\wa=t_\eta 
w$ with $w\in W$ and $\eta\in  Q^\vee$.Then for $i\in \pm[3]$  
$$i^{\wa_*}=\lan\wa(\rho),e_i\ran=\varepsilon_i\lan\eta, 
e_i\ran+\varepsilon_i\lan w(\rho), e_i\ran=\varepsilon_i(\eta, 
e_i)+i^\ww.$$Let also $\ua\in \Wa$, $\ua=t_\tau u$ with $u\in W$ and $\tau\in  
Q^\vee$.Then $$\begin{array}{l}(i^{\wa_*})^{\ua_*}=\varepsilon_i\lan\eta, 
e_i\ran+i^{\ww\ua_*}=\varepsilon_i\lan\eta, 
e_i\ran+\varepsilon_{\iw}\lan\tau,e_{\iw}\ran+i^{\ww 
u_*}\\=\varepsilon_i\lan\eta, e_i\ran+\varepsilon_i\lan\tau,w^{-
1}e_i\ran+i^{\ww u_*}=\varepsilon_i\lan\eta+w(\tau), e_i\ran+i^{\ww 
u_*}.\end{array}$$On the other hand we have $\wa\ua=t_{\eta+w(\tau)}wu$, 
hence  $$i^{(\wa\ua)_*}=\varepsilon_i\lan \eta+w(\tau), 
e_i\ran+i^{(wu)_*},$$and since $i^{(wu)_*}=i^{\ww u_*}$, we finally obtain 
that $i^{(\wa\ua)_*}=i^{\wa_*\ua_*}$. Thus we have that $\Wa$ is isomorphicto 
the subgroup $\Wa_*=\{\ww\mid w\in\Wa\}$ of permutations of $\ganz$.\par
For 
$a\in \ganz$ let $\ov a$ be  the representative of $a\mod 8$ in $[-3,4]$. Then 
using the explicit description of $Q^\vee$ given above, we obtain the following 
permutation representation of $\Wa$.
\begin{theorem}\label{G2}If $\D$ is of 
type $G_2$, then $\Wa$ is isomorphic to the group of all permutations $f$ of 
$\ganz$ such that \item(1) $(-z)^f=-z^f$ for all $z\in \ganz$; \item (2) 
$(z+8k)^f=z^f+8k$ and $(4k)^f=4k$ for all $z, k\in \ganz$ \item (3) $-1^f-
2^f+3^f=0$, $\{\ov {-1^f}, \ov {-2^f}, \ov {3^f}\}=\{-1,-2,3\}$ or $\{\ov {-
1^f}, \ov {-2^f}, \ov {3^f}\}=\{1,2,-3\}$, and $-(1^f-\ov {1^f})\equiv -(2^f-
\ov {2^f})\equiv (3^f-\ov {3^f})\mod 3$.
\end{theorem}
 
\begin{proof}
The statement follows directly from the above discussion.
\end{proof}

\begin{rem}\label{segnogdue}
From the explicit description of $\a_1$, it is clear that ${s_1}_*$ 
acts on $[-3,4]$ as $(1,2)(-1,-2)$.
For $s_2$ we have $s_2(\rho)=\rho-\a_2=e_1-3e_2+2 e_3$, hence ${s_2}_*$
acts on $[-3,4]$ as $(1,-1)(2,3)(-2,-3)$. 
\par
For $w\in W$, let $|\ww|$ be the 
permutation of $[3]$ defined  by $i^{|\ww|}=|i^\ww|$, for $i=1, 2, 3$.
Then from the explicit description of $s_{1*}$ and $s_{2*}$ it is clear that,
for $w\in W$, the parity of $\ell(w_*)$, and hence of $\ell(w)$,  
is exactly the sign of $|\ww|$.  This observation, combined with Lemma 
\ref{lunghezzatau}
and the discussion developed in Subsection \ref{tipoG}, solves the problem 
of determining explicitly $\chi_\l(a),\,\l\in P_{alc}$. With this identification $(-1)^{\ell(v)}$
is the sign of the permutation $|v_*|$ hence, if $\l\in P_{alc}$  and we write
$\l+\r=\mu=\tau+\sum_i r_i e_i$ as described in Subsection \ref{tipoG}, then 
$\chi_\l(a)$ is the sign of the permutation $i\mapsto |r_i|,\,i=1,2,3$.
\par
Finally, we have $s_0(\rho)=\th^\vee+s_\th(\rho)=8\th+\rho-3\th=
\rho+5\th=-6e_1-7e_2+13 e_3$, hence ${s_0}_*$ is the unique permutation $f$ of 
$\ganz$ which has properties (1) and (2) of \ref{G2} and  such that
$1^f=6$, $2^f=7$, $3^f=13$. 
\end{rem}
\vskip10pt
\centerline{\bf ACKNOWLEDGEMENT }
\vskip15pt
We would like to thank H. Eriksson for providing us a copy of his Ph.D. Thesis
\vskip10pt

\providecommand{\bysame}{\leavevmode\hbox to3em{\hrulefill}\thinspace}
\vskip20pt

\footnotesize{
\noindent{\bf 
P.C.}: Dipartimento di Scienze, Universit\`a di Chieti-Pescara, Viale Pindaro 
42, 65127 Pescara,ITALY;\\ {\tt cellini@sci.unich.it}
\par\noindent{\bf P.MF.}: 
Politecnico di Milano, Polo regionale di Como, Via Valleggio 11,22100 Como, 
ITALY;\\ {\tt frajria@mate.polimi.it}
\par\noindent{\bf P.P.}: Dipartimento di 
Matematica, Universit\`a di Roma ``La Sapienza",P.le A. Moro 2, 00185, Roma , 
ITALY;\\{\tt papi@mat.uniroma1.it}}
\end{document}